\documentclass[reqno,12pt,draft]{amsart}
\usepackage[cp1251]{inputenc}
\usepackage[english, russian]{babel}
\usepackage{amsmath,amsthm,amssymb}
\usepackage{multirow}
\begin{document}

\date{18.01.2017}

\author{L.\,M.~Arutyunyan}
\address{Lomonosov MSU, the faculty of mechanics and mathematics, Moscow city.}
\email{Lavrentin@ya.ru}

\title{On growth of the number of determinants with restricted entries}

\markboth{L. M. Arutyunyan}{On growth of the number of determinants}

\maketitle

\begin{small}
\ \ \ \ \ We study the problem about the number of distinct determinants of matrices
with entries from

\ \ \ \ \ a fixed set.

\ \ \ \ \ Bibliography: 10 titles.

\end{small}

\begin{small}
\ \ \ \ \ key words: sum-product phenomenon, determinants.
\end{small}


\section{Introduction}

Let $A$ be a finite subset of a field $\mathbb{F}$ and $D_n(A)$
be a set of all matrices with entries in $A$, namely
$$
D_n(A)=\{D\in \mathbb{F}\ |\ \exists a_{ij}\in A, 1 \le i,j \le n, \det\bigl((a_{ij})\bigr)=D\},
$$
where the symbol $(a_{ij})$ defines the matrix with elements $a_{ij}$.
How big is the size of the set $D_n(A)$ comparing to the size of the set $A$?

The problem we consider is a particular case of the following question
which is quite typical in additive combinatorics.
One considers a function of several variables and explores
how big is the image of the function while the arguments run along a finite set $A$, see \cite{TV}, \cite{61}.

on matrices and distributions 
(particularly, on distributions of their determinants)

Some related problems are considered in papers
\cite{PhV}, \cite{V}, \cite{ShA}, particularly a problem on the distribution of determinants. A continuous counterpart of the examining problem is presented in \cite{GIM}.
Sizes of $|D_n(A)|$ with $n=3,4$ were studied in \cite{C}. For instance, it was proved that
the condition $|A|>\sqrt q$ implies $|D_3(A)|>q/2$, $D_4(A)=\mathbb{F}_q$ (here $q$ is a power of a prime number,
and $\mathbb{F}_q$ is the field of order $q$).
Some other connected questions were also studied there.
Moreover, the set $D_2 (A)=AA-AA$ was intensively studied recently, see \cite{61},\cite{R}
and further references there.

It was proved in paper \cite{PhV} that for an arbitrary $A$ which is a subset of the field $\mathbb{F} = \mathbb{F}_p$,
one has
$$
|D_n(A)| \ge \min(|A|^{3+\frac{1}{45}},p),
$$
There are also some close results.
We prove that for $\mathbb{F}=\mathbb{F}_p$ and arbitrary $A$,
the value of the power unboundedly grows with the size of matrices, more precisely
$$
|D_n(A)| \ge \frac{1}{8} \min(|A|^{c \log n},p),
$$
where $c>0$ is an effective constant. Particularly constant $c=\frac{1}{10}$ is suitable.
The theorem remains true for an arbitrary field of characteristic zero
(of course one can consider $p$ on the right-hand side to be $+\infty$).

The same estimate remains true for the set of permanents instead of the set of determinants, see Remark $3$.

The author expresses thanks to I. D. Shkredov for the formulation of the problem and
a great support throughout the whole research.

\ \

\section{Main definitions.}

For every sets $A, B$, natural numbers $m, n$ and an element of a field $a_0\in \mathbb{F}$,
the following operations are defined:
$$
A+B = \{a+b\ | \ a \in A, b \in B\}, \ \ AB = \{ab\ | \ a \in A, b \in B\}, \ \ a_0*A = \{a_0 a, a\in A\},
$$
$$mA=\{a_1+ a_2 + \ldots+a_n \ | \ a_1,\ldots,a_n \in A\}, \ \
A^m = \{a_1 a_2 \ldots a_n \ | \ a_1,\ldots,a_n \in A\}.$$
The symbol $0_n$ denotes the zero matrix of size $n\times n$.

\section{Proof of the main result}
At the beginning we want to reduce our problem to the case when a set $A$ includes numbers $0$ and $1$.
We need the following lemma to this purpose.

{\bf Lemma 1.} Let $|A|\ge 2$, then $D_{2n}(A) \supset b_0 * D_n(A-A)$ for some $b_0\in \mathbb{F}\setminus\{0\}$.

{\bf Proof.} Let $M_0$ be a matrix $n\times n$ with entries in $A$ such that
$\det(M_0)\neq 0$. As an $M_0$ one can always pick a matrix of the form
\[ \begin{pmatrix}
b & b & b &\ldots &  b & b \\
a & b & a &\ldots & a & a \\
a & a & b &\ldots & a & a \\
\vdots & \vdots & \vdots & \ddots & \vdots & \vdots \\
a & a & a &\ldots & b & a \\
a & a & a &\ldots & a & b
\end{pmatrix}, \]
where $b\in A\setminus\{0\}$, $a\in A\setminus \{b\}$.
Let now $M_1, M_2$ be matrices of size $n$ with entries from 
the set $A$.
Then $D_{2n}$ contains the determinant of the following block matrix: 
\[
\begin{pmatrix}
M_0 & M_1 \\ M_0 & M_2
\end{pmatrix}.
\]
The determinant of this matrix is equal to the determinant of
the difference of the matrices $M_2$ and $M_1$
multiplied by the determinant of the matrix $M_0$.
Indeed, it is easy to get a corner of zeros:
\[ \mathrm{det}\begin{pmatrix}
M_0 & M_1 \\ M_0 & M_2
\end{pmatrix} =
\mathrm{det}\begin{pmatrix}
M_0 & M_1 \\
0_n & M_2-M_1
\end{pmatrix} =\mathrm{det}(M_0) \mathrm{det}(M_2-M_1).
\]
That is why we have the inequality ${D_{2n}(A)\supset \mathrm{det}(M_0)*D_n(A-A)}$.

\ \

{\bf Corollary 1.} Let $|A|\ge 2$. Then there is a set $A'$ with the following properties:
${A'=-A', A'\supset\{0,1\},
|A'|\ge |A|}$,
and also the inclusion
$|D_{2n}(A)|\ge |D_n(A')|$ holds.

{\bf Proof}. As an $A'$ one can consider the set ${(a_0)^{-1} * (A-A)}$ where $a_0$
is an arbitrary element of the set $(A-A)\setminus\{0\}$.
By the previous lemma we have $D_{2n}(A)\supset b_0* D_n(A-A)=b_0 (a_0)^n * D_n (A')$.

\ \

{\bf Theorem 1.} Let $A=-A, A\supset\{0,1\}$. Then we have $D_{m(n-1)+1} (A)\supset n A^m$ for every $m,n\in \mathbb{N}$.

{\bf Proof}. Before the general case of arbitrary $m,n$ we consider $m=3$
and $n=2,3,4$ (it is enough to take a diagonal matrix for $n=1$).
\[ n=2: \ \ \ \ \mathrm{det}\begin{pmatrix}
0 & b_1 & b_2 & b_3 \\
a_1 & 1 & 0  & 0 \\
a_2 & 0 & 1 & 0 \\
a_3 & 0 & 0 & 1
\end{pmatrix} = -(a_1 b_1 + a_2 b_2 + a_3 b_3).
\]
Now we take $a_i,b_i, i=1,2,3$ to be arbitrary elements of $A$, so we get $D_4(A)\supset 3 A^2$.
\[ n=3: \ \ \ \
\mathrm{det}\begin{footnotesize}
\begin{pmatrix}
0 & c_1 & 0 & c_2 & 0 & c_3 & 0 \\
0 & 1 & b_1 & 0 & 0 & 0 & 0 \\
a_1 & 0 & 1 & 0 & 0 & 0 & 0 \\
0 & 0 & 0 & 1 & b_2 & 0 & 0 \\
a_2 & 0 & 0 & 0 & 1 & 0 & 0 \\
0 & 0 & 0 & 0 & 0 & 1 & b_n \\
a_3 & 0 & 0 & 0 & 0 & 0 & 1 \\
\end{pmatrix} = a_1 b_1 c_1 + a_2 b_2 c_2 + a_3 b_3 c_3 \end{footnotesize}
,
\]
Substituting different $a_i,b_i,c_i \in A$ in this formula, we get $D_{3(2-1)+1}(A)\supset 3A^3$.

Now let us consider a matrix for $3 A^4$:
\[
\mathrm{det}\begin{footnotesize}
\begin{pmatrix}
0 & d_1 & 0 & 0 & d_2 & 0 & 0 & d_3 & 0 & 0 \\
0 & 1 & c_1 & 0 & 0 & 0 & 0 & 0 & 0 & 0 \\
0 & 0 & 1 & b_1 & 0 & 0 & 0 & 0 & 0 & 0 \\
a_1 & 0 & 0 & 1 & 0 & 0 & 0 & 0 & 0 & 0 \\
0 & 0 & 0 & 0 & 1 & c_2 & 0 & 0 & 0 & 0 \\
0 & 0 & 0 & 0 & 0 & 1 & b_2 & 0 & 0 & 0 \\
a_2 & 0 & 0 & 0 & 0 & 0 & 1 & 0 & 0 & 0 \\
0 & 0 & 0 & 0 & 0 & 0 & 0 & 1 & c_3 & 0 \\
0 & 0 & 0 & 0 & 0 & 0 & 0 & 0 & 1 & b_3 \\
a_3 & 0 & 0 & 0 & 0 & 0 & 0 & 0 & 0 & 1 \\
\end{pmatrix} \end{footnotesize} =-(a_1 b_1 c_1 d_1+a_2 b_2 c_2 d_2+a_3 b_3 c_3 d_3).
\]

For $m A^n$, one can write down a necessary matrix in the following way.
Let $a=\{a_1,\ldots,a_n\}$ where $a_i=(a_{i,1},\ldots,a_{i,n})$.
Let us define a matrix $M(a_i)$:
\[ M(a_i) =
\begin{pmatrix}
1 & a_{i,n-1} & 0 & \ldots & 0 & 0 \\
0 & 1 & a_{i,n-2} & \ldots & 0 & 0\\
0 & 0 & 1 &         \ldots & 0 & 0\\
\vdots & \vdots & \vdots & \ddots & \vdots & \vdots\\
0 & 0 & 0 & \ldots & 1 & a_{i,2}\\
0 & 0 & 0 & \ldots & 0 & 1\\
\end{pmatrix}.
\]
Now we can define a block matrix $\mathcal{M}(a)$:
\[ \mathcal{M}(a) =
\left(
\begin{tabular}{c|ccccc|c|ccccc} 
0 & $a_{1,n}$ & 0 & 0 & \ldots & 0 & \ldots & $a_{m,n}$ & 0 & 0 & \ldots & 0  \\ \hline
0 & \multicolumn{5}{c|}{\multirow{5}*{\Huge{$M(a_{1})$}}} & \ldots & \multicolumn{5}{c}{\multirow{5}*{\Huge{$0_{n-1}$}}} \\
0 &  &  &  &  &  & \ldots &  &  &  &  &   \\
\vdots &  &  &  &  &  & \ldots &  &  & &  &   \\
0 &  &  &  &  &  & \ldots &  &  &  &  &   \\
$a_{1,1}$ &  &  &  &  &  & \ldots &  &  &  &  &   \\ \hline
\vdots & \vdots & \vdots & \vdots & \vdots & \vdots & $\ddots$ & \vdots & \vdots & \vdots & \vdots & \vdots  \\ \hline
0 &\multicolumn{5}{c|}{\multirow{5}*{\Huge{$0_{n-1}$}}}& \ldots & \multicolumn{5}{c}{\multirow{5}*{\Huge{$M(a_{m})$}}} \\
0 &  &  &  &  &  & \ldots &  &  &  &  &   \\
\vdots &  &  &  &  &  & \ldots &  &  & &  &   \\
0 &  &  &  &  &  & \ldots &  &  &  &  &   \\
$a_{m,1}$ &  &  &  &  &  & \ldots &  &  &  &  &   \\ 
\end{tabular} \right)
\]
Then
$$\det(\mathcal{M}(a)) = (-1)^{n+1}\sum_{i=1}^m \prod_{j=1}^{n} a_{i,j}.$$
Indeed, all non-zero elements of the matrix $\mathcal{M}(a)$ in
the rows from $(i-1) (n-1) + 2$ to $i (n-1)+1$ lie only in the first column and in the block with $M(a_i)$.
Upon that we can choose only one of $m$ elements $a_{i,1}$ from the first row.
Let the element $a_{i_0,1}$ is chosen from the first row.
Then the block with
$M(a_{i_0})$ has not any element in the last row that we can pick,
so we have to pick the element $a_{i_0,2}$ from the last column.
Further there is no element we can pick from the third column in the third right column
in the last two rows, so we indispensably chose $a_{i_0, 3}$. Eventually
we pick the element $a_{i_0, n-1}$ from the second column of the block with $M(a_{i_0})$.
Now in the column of the original matrix with the index $(i_0-1)(n-1) + 1$,
there is the only opportunity to pick $a_{i_0,n}$.
The product of the entries already chosen
equals
$\prod_{j=1}^n a_{i_0,j}$.
Now if one removes the rows and columns with already chosen entries from the matrix,
it remains an upper triangular matrix with all diagonal elements equal to $1$.
It is also easy to assure that the sign of the permutation corresponding to the chosen entries is always the same
(however this is not essential for our further goals as the set $A$ is centrally symmetric).
Thus, if we pick $a$ with every $a_{i,j} \in A$ we get $D_{m(n-1)+1} \supset
m A^n$.

\ \

{\bf Remark 1.} It is easy to see that the following generalization of the previous statement takes place.
In the assumptions of the previous theorem, let
$${k=m_1 (n_1-1)+m_2(n_2-1)+\ldots m_j (n_j-1)},$$
then
$$
D_k(A) \supset m_1 A^{n_1} + m_2 A^{n_2} + \ldots + m_j A^{n_j}.
$$

\ \

{\bf Corollary 2.} For an arbitrary set $A$ we have
$$|D_{2 (m(n-1) + 1)} (A)| \ge |m (A-A)^n|.$$

\ \

Now the main result follows from a result from
The main result will be provided by the following lemma, see \cite{GlK}.

{\bf Lemma 2.} For an arbitrary $A\subset \mathbb{F}_p$, there is an estimate
$$|8^n A^n - 8^n A^n|\ge \frac{1}{8} \min(|A|^n,p).$$
{\bf Proof.} In section $5$ of paper \cite{GlK}, it is proved that
if $|A|\ge 5$ and $N_n=\frac{5}{24}4^n-\frac{1}{3}$ then
$$|N_n A^n-N_n A^n| \ge \frac{3}{8} \min(|A|^n,\frac{p-1}{2}).$$
The announced estimate is obvious for $|A|=0$ and $|A|=1$, while for $|A|\ge 2$, by Cauchy--Davenport Theorem we have $|4A|\ge 5$,
so
$$|8^n A^n - 8^n A^n|\ge |4^n N_n A^n- 4^n A^n|\ge |N_n (4A)^n-N_n (4A)^n|\ge \frac{3}{8} \min(|A|^n,\frac{p-1}{2}),$$
which gives the desired, because $p-1\ge p/3$ for $p\ge 3$ and the case $p=2$ is trivial.

\ \

Lemma 2 together with Corollary 2 provides
$D_{8^{n+1}}(A) \ge \frac{1}{8} \min(|A|^n,p)$, that is why the main result holds true.

{\bf Corollary 3.} $D_n (A) \ge \frac{1}{8} \min(|A|^{0,1 \log n}, p)$

{\bf Proof.} It is easy to see that $|D_n(A)|\ge \frac{1}{8} \min(|A|^{\frac{\log n}{\log 8}-2})$.
For $n \le 2^{10}$ we have $\frac{1}{10} \log n \le 1$, and for $n\ge 2^{10}$ we have
$
\frac{\log n}{\log 8}-2 \ge \frac{1}{10} \log n.
$

\ \

{\bf Remark 2.} The mentioned result remains true for fields of characteristics 0,
since statements like Lemma 2 remains true (moreover, their proofs become easier).

\ \

{\bf Remark 3.} In papers \cite{PhV}, \cite{V1}, a problem similar to our was considered, but for permanents.
In particular, there it was proved that the number of distinct permanents of matrices with entries in a set
$A$ is at least $|A|^{2-\frac{1}{6}+o(1)}$,
where $o(1)$ tends to zero with the growth of matrices size. It is not hard to see that
results achieved here provides the same estimate as in Corollary 3 but for permanents.
Indeed, matrices appeared in Theorem 1 have the same sign of permutations with non-vanishing
product of corresponding elements. So permanent of the matrices appeared there
might differ from their determinants only in the sign.

\ \

The following can be proved analogously to Corollary 3.

{\bf Corollary 4.} Let $\delta \in (0,1)$. Then
$D_n(A)=\mathbb{F}_p$ if $|A|\ge p^{\delta}$ and $n\ge 8 e^{10 \delta}$.

\ \

Obviously, $|D_n(A)|\le |A|^{n^2}$. The following example shows that the upper estimate
can be much stronger than the trivial. For simplicity, let us consider it in
a field of characteristics zero.

{\bf Example 1.} If an estimate of the form $|D_n(A)| \ge C(n) |A|^{n^{\alpha}}$ with some $C(n)>0$
is true for every set $A\subset \mathbb{R}$
then $\alpha$ must not be greater than 1.
Indeed, we can consider $A=\{1,\ldots, m\}$, then as $A^n\subset [1,\ldots, m^n]$,
we have $D_n(A)\subset [-n! m^n, n! m^n]$, so $|D_n(A)|\le C'(n) m^n \le C'(n) |A|^n$.

\end{document}